\documentclass[12pt]{article}
\usepackage{amssymb,amsthm}

\newtheorem*{thm}{Theorem}
\def\card#1{\vert #1 \vert}
\def\gpindex#1#2{\card {#1\colon #2}}

\def\gpcen#1{{\bf Z} (#1)}
\def\ker#1{{\rm ker} (#1)}
\begin{document}
%Notes - Mark L. Lewis - April 2008, revised September 2008.

%These some ideas that came to me while refereeing
%Nenciu's paper: "Brauer pairs of VZ-groups."

\title{Brauer pairs of Camina $p$-groups of nilpotence class $2$}

\author {
       Mark L.\ Lewis
    \\ {\it Department of Mathematical Sciences, Kent State University}
    \\ {\it Kent, Ohio 44242}
    \\ E-mail: lewis@math.kent.edu
       }
\date{September 26, 2008}

\maketitle

\begin{abstract}
In this paper, we find a condition that characterizes when two
Camina $p$-groups of nilpotence class $2$ form a Brauer pair.

MSC primary: 20C15
\end{abstract}

\section{Introduction }

Throughout this note, all groups are finite.  Two nonisomorphic
groups $G$ and $H$ are said to form a {\it Brauer pair} if $G$ and
$H$ have identical character tables and identical power maps.  In
\cite{Brauer}, Brauer had asked if there exist any such pairs. The
first examples of Brauer pairs were found by Dade in \cite{Dade}.
Other examples of Brauer pairs can be found in the \cite{EiMu} and
\cite{NenBr}.  In this paper, we find a condition that characterizes
when Camina groups of nilpotence class $2$ form Brauer pairs.

In \cite{VZgroups}, we defined a {\it VZ-group} to be a group where
all the nonlinear irreducible characters vanish off the center, and
in that paper, we characterizes when two such groups have identical
character tables.  Using our characterization of the character
tables of VZ-groups, A. Nenciu was able to characterize those
VZ-groups that were Brauer pairs (see \cite{NenVZ}).

A group $G$ is a Camina group if for all $g \in G \setminus G'$ we
have that $g$ has $gG'$ as its conjugacy class.  The nilpotent
Camina groups of nilpotence class $2$ are VZ-groups.  These have
also been studied under the name semi-extraspecial groups in
\cite{Beis} and \cite{Ver}.  (One can see in \cite{extreme} and
\cite{Norit} that the condition semi-extraspecial is equivalent to
the condition of being a Camina $p$-group of nilpotence class $2$.)
Our characterization in \cite{VZgroups}
%of VZ-groups with identical character tables
simplified for Camina groups of nilpotence class $2$ to say: let $G$
and $H$ be nilpotent Camina groups of nilpotence class $2$. Then $G$
and $H$ have identical character tables if and only if $\gpindex
G{G'} = \gpindex H{H'}$ and $\card {G'} = \card {H'}$.  Our goal in
this note is to show that Nenciu's characterization of Brauer pairs
of VZ-groups also has an easy simplification to Camina groups of
nilpotence class $2$.

Any nilpotent Camina group must be a $p$-group for some prime $p$,
and any VZ-group will be the direct product of a $p$-group and an
abelian $p'$-group. Also, Nenciu showed in \cite{NenVZ}, that there
are no Brauer pairs of VZ-groups that are also $2$-groups.
Therefore, we may assume that we are working with Camina $p$-groups
where $p$ is an odd prime. Finally, if $P$ is a $p$-group, then we
have the subgroups $\mho_1 (P) = \langle x^p \mid x \in P \rangle$
and $\Omega_1 (P) = \langle x \in P \mid x^p = 1 \rangle$. When $p$
is odd and $P$ has nilpotence class $2$, it is known that in fact
$\mho_1 (P) = \{ x^p \mid x \in P \}$.

With $\mho_1 (P)$ in hand, we can state our simplified condition
that characterizes Brauer pairs of Camina $p$-groups with nilpotence
class $2$.  It shows that if $P$ and $Q$ are Camina $p$-groups of
nilpotence class $2$, then $P$ and $Q$ form a Brauer pair if and
only if $\gpindex P{P'} = \gpindex Q{Q'}$, $\card {P'} = \card
{Q'}$, and $\card {\mho_1 (P)} = \card {\mho_1 (Q)}$.  In other
words, the only hypothesis needed to have a Brauer pair beyond the
conditions to have identical character tables is that $\mho_1 (P)$
and $\mho_1 (Q)$ have the same order.

\begin{thm}[Main Theorem]
Let $P$ and $Q$ be nonisomorphic Camina $p$-groups of nilpotence
class $2$ for some odd prime $p$.  Then $P$ and $Q$ form a Brauer
pair if and only if $\gpindex P{P'} = \gpindex Q{Q'}$, $\card {P'} =
\card {Q'}$, and $\card {\mho_1 (P)} = \card {\mho_1 (Q)}$.
\end{thm}

We conclude this introduction by noting that we will show the
condition that $\card {\mho_1 (P)} = \card {\mho_1 (Q)}$ is
equivalent to $\gpindex P{\Omega_1 (P)} = \gpindex Q{\Omega_1 (Q)}$
for these groups.  Hence, the Main Theorem could also be restated as
saying that when $P$ and $Q$ are nonisomorphic Camina $p$-groups of
nilpotence class $2$, then $P$ and $Q$ form a Brauer pair if and
only if $\gpindex P{P'} = \gpindex Q{Q'}$, $\card {P'} = \card
{Q'}$, and $\gpindex P{\Omega_1 (P)} = \gpindex Q{\Omega_1 (Q)}$

\section{Results}

%I need to add exposition here.
Before we prove the main theorem, we need to state some notation.
Let $P$ be a $p$-group with nilpotence class $2$, where $p$ is an
odd prime. We define $\nu_P : P \rightarrow P$ by $\nu_P (x) = x^p$.
When $P'$ is elementary abelian, it is not difficult to see that
$\nu_P$ is a homomorphism of $P$ and its image is $\mho_1 (P)$.  The
fact that $p$ is odd and $P'$ is central and elementary abelian are
necessary for $\nu_P$ to be a homomorphism.  Also, $\ker {\nu_P} =
\{ x \in P \mid x^p = 1 \} = \Omega_1 (P)$.  Hence, we have
$\gpindex G{\Omega_1 (P)} = \card {\mho_1 (P)}$.  As $P'$ is
elementary abelian, $P' \le \ker {\nu_P}$, and if $P/\gpcen P$ is
elementary abelian, then $\mho_1 (P) \le \gpcen P$. If $P$ is a
VZ-group, both of these occur; so when $P$ is a VZ-group, we can
view $\nu_P : P/P' \rightarrow \gpcen P$ as a homomorphism.  We can
also define the projection map $\phi_P : \gpcen P \rightarrow P/P'$
by $\phi_P (z) = zP'$ for every $z \in \gpcen P$.

In \cite{NenVZ}, Nenciu showed that if $P$ and $Q$ are nonisomorphic
VZ-groups associated with the prime $p$, then $P$ and $Q$ form a
Brauer pair if and only if there exist isomorphisms $\hat\alpha :
P/P' \rightarrow Q/Q'$ and $\hat\beta : \gpcen P \rightarrow \gpcen
Q$ so that $\hat\alpha \circ \phi_P = \phi_Q \circ \hat\beta$ and
$\nu_Q \circ \hat\alpha = \hat\beta \circ \nu_P$.

\begin{proof}[Proof of Main Theorem]
We know when $P$ is Camina group with nilpotence class $2$ that
$\gpcen P = P'$.  It follows that $\phi_P$ the trivial map.
Similarly, $\phi_Q$ will be the trivial map.  Thus, Nenciu's
condition for $P$ and $Q$ to be a Brauer pair is equivalent to
finding isomorphisms $\hat\alpha$ and $\hat\beta$ so that $\nu_P
\circ \hat\alpha = \hat\beta \circ \nu_Q$.

Suppose that $P$ and $Q$ form a Brauer pair.  Then isomorphisms
$\hat\alpha$ and $\hat\beta$ exist. We already know that this
implies that $\gpindex P{P'} = \gpindex Q{Q'}$ and $\card {P'} =
\card {Q'}$.  We have $\hat\beta \circ \nu_P (P/P') = \hat\beta
(\nu_P (P/P')) = \hat\beta (\mho_1 (P))$, and $\nu_Q \circ
\hat\alpha (P/P') = \nu_Q (\hat\alpha (P/P')) = \nu_Q (Q/Q') =
\mho_1 (Q)$.  It follows that $\hat\beta (\mho_1 (P)) = \mho_1 (Q)$,
and we conclude that $\card {\mho_1 (P)} = \card {\mho_1 (Q)}$ as
desired.

We now suppose that $\gpindex P{P'} = \gpindex Q{Q'}$, $\card {P'} =
\card {Q'}$, and $\card {\mho_1 (P)} = \card {\mho_1 (Q)}$. Now,
$\mho_1 (P) \le \gpcen P = P'$, so $\mho_1 (P)$ is an elementary
abelian $p$-group.  Similarly, $\mho_1 (Q) \le Q'$ is an elementary
abelian $p$-group.  Since $\mho_1 (P)$ and $\mho_1 (Q)$ have the
same size, they are isomorphic.  Write $\beta$ for the isomorphism
from $\mho_1 (P)$ to $\mho_1 (Q)$.  Since $P'$ and $Q'$ are
elementary abelian of the same size, we can extend $\beta$ to an
isomorphism $\hat\beta$ from $P'$ to $Q'$.

We know that $P' \le \Omega_1 (P)$ and $Q' \le \Omega_1 (Q)$. Since
$P/P'$ is an elementary abelian $p$-group, there is a subgroup $A$
of $P$ so that $P/P' = \Omega_1 (P)/P' \times A/P'$.  Similarly,
there is a subgroup $B$ of $Q$ so that $Q/Q' = \Omega_1 (Q)/Q'
\times B/Q'$. Now, $\Omega_1 (P)$ is the kernel of $\nu_P$, so we
know that $\nu_P : P/\Omega_1 (P) \cong \mho_1 (P)$. Restricting
$\nu_P$ to $A$, we have that $\nu_P : A/P' \cong \mho_1 (P)$.
Similarly, $\nu_Q : B/Q' \cong \mho_1 (Q)$. We can define $\alpha :
A/P' \rightarrow B/P'$ to be the unique map so that $\nu_P \circ
\alpha = \beta \circ \nu_Q$. Observe that $\Omega_1 (P)/P'$ and
$\Omega_1 (Q)/Q'$ are elementary abelian groups of the same size.
Hence, there is an isomorphism $a : \Omega_1 (P)/P' \rightarrow
\Omega_1 (Q)/Q'$.  Note that $\nu_P \circ a = 1$.  We can define an
isomorphism $\hat\alpha : P/P' \rightarrow Q/Q'$ by $\hat\alpha =
(\alpha,a)$.  We observe that $\nu_P \circ \hat\alpha = \nu_P \circ
\alpha = \beta \circ \nu_Q = \hat\beta \circ \nu_Q$.
\end{proof}

Now, suppose that $p$ is an odd prime, and let $P$ be a Camina
$p$-group with nilpotence class $2$ and $\card {P'} = p^n$.  Thus,
there are $n+1$ choices for the value of $\card {\mho_1 (P)}$.  If
the number of Camina groups with $\card {P'} = p^n$ and $\gpindex
P{P'} = p^{2m}$ with $m \ge n$ is bigger than $n+1$, then this will
be another source of Brauer pairs.

We will have that $P$ and $Q$ form a Brauer pair when $P$ and $Q$
are Camina groups of nilpotence class $2$ and exponent $p$ with
$\gpindex P{P'} = \gpindex Q{Q'}$ and $\card {P'} = \card {Q'}$.
Verardi has shown in \cite{Ver} a number of ways of constructing
nonisomorphic Camina $p$-groups with nilpotence class $2$ and
exponent $p$ for a fixed $\gpindex P{P'}$ and $\card {P'}$.  Thus,
we get a number of Brauer pairs this way.

If $n = 1$, then $P$ is extra-special of order $p^{2m+1}$.  We know
that there are $2$ such groups one with exponent $p$ (and hence,
$\card {\mho_1 (P)} = 1$) and one with exponent $p^2$ (and hence,
$\card {\mho_1 (P)} = p$).  Thus, we do not get any Brauer pairs
here.  Of course, this fact is well-known.

Next, we decided to look at the Camina groups with $\gpindex P{P'} =
p^4$ and $\card {P'} = p^2$ using MAGMA and the library of small
groups.  Notice that there are $3$ possibility for the order of
$\mho_1 (P)$.  We have checked the primes up to $31$. For each
prime, we have found that there are $p+3$ Camina groups with
$\gpindex P{P'} = p^4$ and $\card {P'} = p^2$.  For each prime that
we tested, one group has exponent $p$, and so, $\card {\mho_1 (P)} =
1$; one group has $\card {\mho_1 (P)} = p$; and the remaining $p+1$
groups have $\card {\mho_1 (P)} = p^2$.  We note that one of the
groups in the last class has that $\Omega_1 (P)$ is not abelian, and
the remaining $p$-groups in that class have $\Omega_1 (P)$ is
abelian.  Thus, we obtain Brauer pairs by taking any two
nonisomorphic Camina groups with nilpotence class $2$ where
$\gpindex P{P'} = p^4$, $\card {P'} = p^2$, and $\card {\mho_1 (P)}
= p^2$.

We note that there are two examples in \cite{NenVZ} which motivated
our study. The examples there actually meet the hypotheses of our
theorem. Both of the examples there have $\card {\mho_1 (P)} = p^2$.
It follows that $\card {\Omega_1 (P)} = p^4$; one of the examples
has $\Omega_1 (P)$ abelian and the other example has $\Omega_1 (P)$
nonabelian.

\end{document}